\documentclass{jnmp01b}

\usepackage{amsmath}

\numberwithin{equation}{section}

\newcommand{\nz}{\mbox{${\Bbb Z}$}}
\newcommand{\cC}{\mbox{${\cal C}$}}
\newcommand{\vphi}{\varphi}
\theoremstyle{plain}
 \newtheorem{thm}{Theorem}[section]
 \newtheorem{prop}[thm]{Proposition}
 \newtheorem{lemma}[thm]{Lemma}
\theoremstyle{definition}
 \newtheorem{defn}{Definition}[section]

\makeatletter
\DeclareRobustCommand{\primfrac}[1]{%
  \PackageWarning{amsmath}{%
Foreign command \@backslashchar#1; %
\protect\frac\space or \protect\genfrac\space should be used instead%
  }
  \global\@xp\let\csname#1\@xp\endcsname\csname @@#1\endcsname
  \csname#1\endcsname
}
\makeatother

\begin{document}

\renewcommand{\evenhead}{J.\ Ding and B.\ Feigin}
\renewcommand{\oddhead}{Quantized W-algebra of ${\frak sl}(2,1) $
and Quantum Parafermions of 
$U_q(\hat {\frak sl}(2))$}


\thispagestyle{empty}

\begin{flushleft}
\footnotesize \sf
Journal of Nonlinear Mathematical Physics \qquad 2000, V.7, N~2,
\pageref{firstpage}--\pageref{lastpage}.
\hfill {\sc Article}
\end{flushleft}

\vspace{-5mm}

\copyrightnote{2000}{J.\ Ding and B.\ Feigin}

\Name{Quantized W-algebra of ${\frak sl}(2,1) $
and Quantum Parafermions of 
$U_q(\hat {\frak sl}(2))$ }

\label{firstpage}

\Author{Jintai DING~$^\dag$ and Boris FEIGIN~$^\ddag$}

\Adress{$^\dag$ Department of Mathematical Sciences, 
        University of Cincinnati, USA\\[2mm]
        $^\ddag$ Landau Institute, Moscow, Russia}

\Date{Received December 15, 1999; Accepted January 28, 2000}

\begin{abstract}
\noindent
In this  paper,  we establish the connection between the 
quantized W-algebra of ${\frak sl}(2,1)$ 
and quantum parafermions of $U_q(\hat {\frak sl}(2))$
that a shifted product of the two quantum parafermions
of $U_q(\hat {\frak sl}(2))$ generates the
quantized W-algebra of ${\frak sl}(2,1)$. 
\end{abstract}


\section{ Introduction} 

The Lie algebra $\hat {\frak sl}(2)$ has a current realization 
that is given by three current operators $e(z)$, 
$h(z)$ and $f(z)$. 
For $\hat {\frak sl}(2)$, when we consider the case of bosonization of 
the current operator, the operators $e(z)$ and $f(z)$ 
can be decomposed in the following way: 
\begin{gather*}
e(z)=v^+(z)\Psi^+(z), \\
f(z)= v^-(z)\Psi^-(z), 
\end{gather*}
where  $\Psi^{\pm}(z)$ are operators that commute with $h(z)$ and 
$v^{\pm}(z)$ are vertex operators generated by the Heisenberg algebra of 
$h(z)$.  
These two operators $\Psi^{\pm}(z)$ are called parafermions. 

 For the quantum affine algebra 
 $U_q(\hat {\frak sl}(2))$, which is a q-deformation of the 
universal enveloping algebra of $\hat {\frak sl}(2)$, Drinfeld 
presented a loop realization  of affine quantum groups with 
current generators 
\cite{1}. This, for   $U_q(\hat {\frak sl}(2))$,   gives us the 
 quantized current operators corresponding to $e(z)$, $h(z)$ and $f(z)$ of 
 $\hat {\frak sl}(2)$, which  are 
$X^+(z)$,  $\vphi(z)$,  $\psi(z)$ and $ X^-(z)$. 
Here $\vphi(z)$ and   $\psi(z)$ correspond to the positive and negative half of 
$h(z)$ respectively. 
For $U_q(\hat {\frak sl}(2))$, similarly 
when we consider the case of bosonization of 
the current operator, the operators $X^+(z)$ and $X^-(z)$ 
can be decomposed in the following way: 
\begin{gather*}
X^+(z)=V^+(z)\Phi^+(z),  \\
X^-(z)= V^-(z)\Phi^-(z), 
\end{gather*}
where  $\Phi^{\pm}(z)$ are operators that commute with 
 $\vphi(z)$,  $\psi(z)$ and 
$V^{\pm}(z)$ are vertex operators generated by the Heisenberg algebra of 
  $\vphi(z)$,  $\psi(z)$\cite{2}\cite{3}\cite{4}.

For a given affine Lie (super) algebra $\hat {\frak g}$, there exist different  ways to 
construct the associated   W-algebra, which is an
extended conformal algebra.  
  One is the well-known  Sugawara construction of Casimir type. 
Another is the Drinfeld Sokolov reduction.  
There is also another type of construction coming  from the description of  the associated W-algebra as 
the algebra 
generated by the current operators which commutes with the screening operators associated to  
  $\hat {\frak g}$. This  is evident for the case $g={\frak sl}(2)$, where the construction  
produces  the Virasoro algebra. 
The key idea is to consider the screening operators associated to  the root of 
$ {\frak g}$ in terms of bosonization, then we  try to find the simplest current operators that 
commute with those screen operators also in terms of bosonization. It turns out 
that it also gives a W-algebra.  
Such a construction was first quantized for the case of Virasoro, which gives us 
the quantized Virasoro algebra\cite{13}
For the case of 
$g={\frak sl}(N)$, such a construction   is also utilized in 
the case of constructing quantized W-algebra  associated to ${\frak sl}(N)$\cite{5}. 
Recently, following   the idea of the 
description of quantum W-algebra using  the quantization of screening operators \cite{5},   
we derived  the quantized W-algebra of ${\frak sl}(2,1)$\cite{6}. 
The generator of 
this algebra $l_1(z)$ is bosonized as a sum of three vertex operators, 
while the screening  operators are a set of  two fermions or a set of one quantized screening operators of 
${\frak sl}(2)$ and a fermion. For these two different sets of 
screen operators,  they actually produce  equivalent $l_1(z)$\cite{6}, where  $l(z)$ is given as  
 an uniquely determined operator that commutes with the two screening operators  up-to a total 
difference.  
We see clearly \cite{6}, in terms of the  description of one quantized screening operators of 
${\frak sl}(2)$ and a fermion, 
the degeneration of the two screen operators, when q goes to 1, are exactly 
the screening operators associated to the Lie super algebra  $\hat {\frak sl}(2,1)$. 
This is the reason the algebra generated by $l_1(z)$ is called the quantized W-algebra of ${\frak sl}(2,1)$.  
 The degeneration of 
this operator $l_1(z)$, when q goes to 1,  still commutes with 
the screening operators, and gives  us  the 
associated  W-algebra of ${\frak sl}(2,1)$. 
 
It is known in mathematics as a folklore that there is connection between the 
W-algebra corresponding to  the Lie super algebra $\hat {\frak sl}(2,1)$ and these two parafermions. Namely that 
in the operator product expansion of $\Psi^{-}(w)\Psi^{+}(z)$, one 
of the components gives us a realization of the  W-algebra of ${\frak sl}(2,1)$, namely 
if we look at the bosonized $\Psi^{-}(w)\Psi^{+}(z)$'s components in terms of product expansion, 
the corresponding component commutes with two fermions, which is 
equivalent to the screening operators of  ${\frak sl}(2,1)$, 
 which are equivalent as shown in \cite{6} for the quantized case. 

In this paper, we will use the bosonization formulas in 
\cite{2} to establish again 
the connection between the quantized W-algebra of 
${\frak sl}(2,1)$ and the two quantum parafermions. 
We show that the operator $l(z)$ can be identified with 
the operator $\int_C \Phi^-(w)\Phi^+(z)dw$, where the $k$ is 
the central element or say level and the contour C is around the 
point $w=zq^{k+2}$. 
These extend  the classical  theory respectively. 

 In Section 2, we will 
introduce  the quantized W-algebra of ${\frak sl(2,1 )} $ and 
some basic facts. 
In Section 3 of this paper, we will introduce the 
basic theory about the algebra $U_q(\hat {\frak sl}(2))$ and 
baosonization and  quantum parafermions.
 In Section 4, we will show the connection
between the  quantized W-algebra of ${\frak sl(2,1 )} $ and 
the quantum parafermions of  $U_q(\hat {\frak sl}(2))$, which is 
obtained as a corollary of our previous work \cite{6}. 

\newpage

\section { Quantized  W-algebra of   ${\frak sl}(2,1)$}

In this section, we introduce the definition of 
the quantized  W-algebra of   ${\frak sl}(2,1)$ 
presented in \cite{6}. This is mainly from the last section of 
\cite{6}.

\begin{defn}
The
Heisenberg  algebra 
$H_{q,p}(2)$ is an associative algebra with generators $a_i[n],
n\in{\Bbb Z}$, and the commutation relations of the generators are: 
\begin{equation*}    
\left[ a_i[n],a_j[m] \right] = \frac{1}{n}A_{ij}(n)
\delta_{n,-m},
\end{equation*}
which is defined on the field of the rational functions of 
 $p$ and $q$, two generic parameters with   $p,q \in {\Bbb C}^*, |q|<1$,
 $A_{i,i}=1$ and $A_{ij}(n)=a(n)$ is
 a rational function of $p$ and $q$ for $i\neq j$. 
\end{defn}

Let  $\mu$  be an  element of a two dimensional 
space  $A_2$ generated by $\alpha_i$, for $i=1,2$. Let 
$\alpha_i^*$ be the generator of its dual space $A_2^*$, such that 
$\alpha_i^*(\alpha_i)=2$ and $\alpha_i^*(\alpha_j)=-b/\beta_i$.
Let $\pi_\mu$ be
the Fock representation of $H_{q,p}(2)$ generated by a vector
$v_\mu ( \mu \in A_2^*)$, such that $a_i[n] v_\mu = 0, n>0$, and $a_1[0] v_\mu =
\mu(\alpha_i) v_\mu$. We assume here $\beta_1$ and $\beta_2$ generic. 

Introduce operators $Q_i$, i=1,2, which satisfy commutation
relations $[a_i[n],Q_i] = 2\beta_i\delta_{n,0}$, 
$[a_j[n],Q_i]=- b\delta_{n,0}$. The operators
$e^{Q_i}$ act from $\pi_\mu$ to $\pi_{\mu+\beta\alpha_i}$.

We define two quantized screening currents as the generating function
\begin{equation*}    
S^+_i(z) = e^{Q_i} z^{s^+_i[0]} :\exp \left( \sum_{m\neq 0} s^+_i(m)
a_i[m]z^{-m} \right):, 
\end{equation*}
for $i=1,2$,  where $s^+_i[m]$ are in the polynomial ring of 
$p,q$ over ${\Bbb C}^*$  
for $m\neq 0$ and $s^+_i[0] = a_i[0]$ and by  $: :$  we mean the normal 
ordered product expansion. 

We  impose the following  equalities as an assumption on these two current operators. 

\medskip
{\bf Assumption:}
\begin{gather*}
S^+_i(z)S^+_i(w)=(z-w):S^+_i(z)S^+_i(w):,  \\[1em]
S^+_2(z)S^+_1(w)=z^{-b}f_{2,1}(z,w):S^+_2(z)S^+_1(w): \\
\qquad = z^{-b}
\exp\left (\exp 1/n \Sigma_{m>0} A_{2,1}s^+_2(m)s^+_i(-m)w^m/z^m)\right )
:S^+_2(z)S^+_1(w):, \\[1em]
S^+_1(z)S^+_2(w)=z^{-b}f_{1,2}(z,w):S^+_1(z)S^+_2(w): \\
\qquad =
z^{-b} \exp\left (\exp \Sigma_{m>0}
A_{1,2}(m) s^+_2(m)s^+_2(-m)w^m/z^m)\right )
:S^+_1(z)S^+_2(w): ,
\end{gather*}
and 
\[\lim_{q\rightarrow 1}  f_{(1,2)}(z,w)= 
 \lim_{q\rightarrow 1}f_{(2,1)}(z,w)=  
(1-w/z)^b.\]
As defined in Definition 2.1, $q$ is a complex number, which allows us to take the limit above.

Clearly 
$S^+_1(z)$ and $S^+_2(z)$ are two fermions, which are not the 
inverse of each other. 

 Let 
$l_1(z)$ be  a current operator in the form:
\[l_1(z)= \Lambda_1(z)+\Lambda_2(z)+\Lambda_3(z),\]  where 
$\Lambda_i(z)$ are  the generating functions: 
\begin{equation*}    
\Lambda_i(z) =g_i
 :\exp \left( \sum_{} \lambda_{ij}(m)
a_j[m]z^{-m} \right):, 
\end{equation*}
 $\lambda_i[m]$ are in ${\Bbb C}[p,q]$  for $i=1,2,3$, $g_1=1$,  
the integration contour is around 0;    
and  $l_1(z)$  commutes  with the 
action of the quantized screening  operators  $\int S^+_i(z)dz/z$. 

Here we can actually also  define $l_1(z)$ even by a 
stronger  condition that it  commutes with the operators 
$S^+_i(z)$ up-to a total difference. 

The main goal is to find out  if such an operator actually exists;  and 
if it exists, we would like to find out if  it is  unique.  

In this case, the simplest situation that we have is 
to impose the following assumptions on the correlation functions. 

\medskip
{\bf Assumption:} 
The  correlation 
functions   between $S_i^+(z)$ and 
$\Lambda_j(w)$ are  $1$, for the two pairs 
$i=1, j=3$, and $i=2,j=1$,  which also means that 
for either  pair of the  operators, they commute with each other.
The correlation functions between $S^+_1(z)$ and $\Lambda_1(z)$
satisfy  the condition that 
the two 
products  $\Lambda_1(z) S^+_1(w)$
and  $S^+_1(w) \Lambda_1(z)$ have the same correlation functions 
and 
\begin{gather*}
\Lambda_1(z) S^+_1(w) =
A\frac{(z-w)}{(z-wpq^{-1})} :\Lambda_1(z)
S^+_1(w):, \quad \quad |z|\gg|w|, \\
S^+_1(w) \Lambda_1(z) =
A\frac{(z-w)}{(z-wpq^{-1})} :\Lambda_1(z)
S^+_1(w):, \quad \quad |w|\gg|z|,
\end{gather*}
and 
\[A=pq^{-1}.\]

\begin{prop}[\cite{6}]
The correlation functions of the products 
$\Lambda_1(z) S^+_1(w)$ and \linebreak[4]
$S^+_1(w) \Lambda_1(z)$ must be equal and the correlation functions 
must have only one pole and one zero. For some number $A'$, 
\begin{gather*}
\Lambda_2(z) S^+_1(w) =
A'\frac{(z-wp'_1)}{(z-wp'_2)} :\Lambda_2(z)
S^+_1(w):, \quad \quad |z|\gg|w|,\\
S^+_1(w) \Lambda_2(z) =
A'\frac{(z-wp'_1)}{(z-wp'_2)} :\Lambda_2(z)
S^+_1(w):, \quad \quad |w|\gg|z|, \\
A'p'_1/p'_2=1. \\
A(1-p_1/p_2)p_2: \Lambda_1(z)
S^+_1(zp_2^{-1}):=- A'(1-p'_1/p'_2)p'_2: \Lambda_2(z)
S^+_1(z{p'_2}^{-1}):. 
\end{gather*}
Let $p=p'_2$ and $q=p'_2/p_2$, then $p'_1=1$,
$g_2=pq^{-1}(pq^{-1}-1) .$
\end{prop} 

Let  
\begin{gather*}
S^+_i(z)\Lambda_j(w)=S\Lambda_{ij}(z,w):S^+_i(z)\Lambda(w):, \\
\Lambda_i(z)S^+_j(w)=\Lambda S_{ij}(z,w):S^+_i(z)\Lambda(w): .
\end{gather*}

We also  impose the  following {\bf Assumption}: for some number $B$, 
\begin{gather*}
\Lambda_2(z) S^+_2(w) =
B\frac{(z-wq_1)}{(z-wq_2)} :\Lambda_2(z)
S^+_2(w):, \quad \quad |z|\gg|w|,\\ 
S^+_2(w) \Lambda_2(z) =
B\frac{(z-wq_1)}{(z-wq_2)} :\Lambda_2(z)
S^+_2(w):, \quad \quad |w|\gg|z|.
\end{gather*}

Then, we have 
\begin{prop}
The correlation functions of the products 
$\Lambda_3(z) S^+_2(w)$
and  $S^+_2(w) \Lambda_3(z)$ must be equal and the correlation functions 
must have only one pole and one zero. For some number $B'$, 
\begin{gather*}
\Lambda_3(z) S^+_2(w) =
B'\frac{(z-wq_1)}{(z-wq'_2)} :\Lambda_3(z)
S^+_2(w):, \quad \quad |z|\gg|w|,\\
S^+_2(w) \Lambda_3(z) =
B'\frac{(z-wq_1)}{(z-wq'_2)} :\Lambda_3(z)
S^+_2(w):, \quad \quad |w|\gg|z|,\\
B'(1-q_1/q_2)q_2: \Lambda_2(z)
S^+_2(zq_2^{-1}):=- B'(1-q_1/q'_2)q'_2: \Lambda_3(z)
S^+_2(z{q'_2}^{-1}):, \\
B'q'_1/q'_2=1.
\end{gather*}
\end{prop}

Let $p'=q_2'/q_2$.

\begin{thm}[\cite{6}]
The  operator $l_1(z)$ exists and is  uniquely determined  if and only if 
\begin{gather*}
q'=q\\
q_2=q_1p, \\
\Lambda S_{2,1}(z, w)
=A \frac {f_{2,1}(z{q'_2}^{-1}, w)} {f_{2,1}(z{q_2}^{-1},w)}=
\frac {f_{1,2}(w, z{q'_2}^{-1})}{f_{1,2}(w, z{q_2}^{-1})}, \\
A^{-1}\Lambda S_{2,1}(wq'_2,z)=S\Lambda_{2,2}(zp^{-1},w) . 
\end{gather*}
\end{thm}

\begin{prop} If $l_1(z)$ exists, then 
\begin{gather*}
{f_{2,1}(z, w)}= \frac 
{(w/z| {q'_2}^{-1}pq, q)_\infty}
{(w/z| {q'_2}^{-1}q^{}, q)_\infty},\\
 {f_{1,2}(w, z)}= 
\frac 
{(z/w| {q'_2}p^{-1}, q)_\infty}
{(z/w| {q'_2}, q)_\infty}.
\end{gather*}
\end{prop}

As explained in \cite{6}, we know that the number $b$ actually does not 
affect the commutation relation of $l(z)$ with itself at all. So this 
family of operators related to the parameter $b$ are actually just 
a kind of rescaling.
We call the associative algebra generated by the Fourier components of   the operator $l_1(z)$, 
{\bf the quantized W-algebra of ${\frak sl}(2,1)$}.

\section{  $U_q({\frak sl}(2))$, its bosonization and
quantum  parafermions}

For the definition, we 
will use directly 
the current realization of $U_q(\hat {\frak sl}(2))$ given by 
\linebreak[4]%
Drinfeld\cite{1}.

\begin{defn}
The algebra $U_q(\hat {\frak sl}(2))$ is an associative algebra with unit 
1 and the generators: $\vphi(m)$,$\psi(-m)$, $X^{\pm}(l)$, for 
$i=1,...,n-1$, $l\in \nz $ and $m\in \nz_{\leq 0}$ and a central
 element $c$. Let $z$ be a formal variable and 
 $X^{\pm}(z)=\sum_{l\in \nz}X^{\pm}(l)z^{-l}$, 
$\vphi(z)=\sum_{m\in \nz_{\leq 0}}\vphi(m)z^{-m}$ and 
$\psi(z)=\sum_{m\in \nz_{\geq 0}}\psi(m)z^{-m}$. In terms of the 
formal variables, 
the defining relations are 
\begin{gather*}
\vphi(0)\psi(0)=\psi(0)\vphi(0)=1, \\
\vphi(z)\vphi(w)=\vphi(w)\vphi(z), \\
\psi(z)\psi(w)=\psi(w)\psi(z), \\
\vphi(z)\psi(w)\vphi(z)^{-1}\psi(w)^{-1}=
  \frac{g(z/wq^{-c})}{g(z/wq^{c})}, \\
\vphi(z)X^{\pm}(w)\vphi(z)^{-1}=
  g(z/wq^{\mp \frac{1}{2}c})^{\pm1}X^{\pm}(w), \\
\psi(z)X^{\pm}(w)\psi(z)^{-1}=
  g(w/zq^{\mp \frac{1}{2}c})^{\mp1}X^{\pm}(w), \\
[X^+(z),X^-(w)]=\frac{1}{q-q^{-1}}
  \left\{ \delta(z/wq^{-c})\psi(wq^{\frac{1}{2}c})-
          \delta(z/wq^{c})\vphi(zq^{\frac{1}{2}c}) \right\}, \\
(z-q^{\pm a}w)X^{\pm}(z)X^{\pm}(w)=
  (q^{\pm a}z-w)X^{\pm}(w)X^{\pm}(z), 
\end{gather*}
where
\[ \delta(z)=\sum_{k\in \nz}z^k, \quad
   g(z)=\frac{q^{a}z-1}{z-q^{a}}\quad \text{(expanded around
 $z=0$) },
\quad a=2. \]
\end{defn}

For this current realization, Drinfeld also gave the Hopf algebra 
structure \cite{7}. 

We first  introduce the bosonic
 representation of $U_q(\hat{\frak sl}(2))$ given in \cite{2}. 
We will use the notation from \cite{2}. 

Let $ \{\alpha_n,\bar \alpha_n,\beta_n | n \in \nz 
\}$ be a set of operators satisfying the
following commutation relations:
\begin {gather*}
[\alpha_m,\alpha_{-m}]=\frac {[2m][km]}  {2},\\
[\bar \alpha_m,\bar \alpha_{-m}]=-\frac {[2m][km]}  {2},\\
[\beta_m,\beta_{-m}]=\frac {[2m][(k +2)m]}  {2}. 
\end{gather*}
The other commutators are zero.
These operators  form the direct sum of three Heisenberg algebras.

We define
\begin{align*}
N_ +&=\cC[\alpha_m,\bar \alpha_m,\beta_m]_{m>0}, \\
N_-&=\cC[\alpha_m,\bar \alpha_m,\beta_m]_{m<0}.
\end{align*}
The left Fock module $F_{l,m_1,m_2}$ is uniquely characterized by the following
properties:
there exists a vector $|l,m_1,m_2>$ in $F_{l,m_1,m_2}$ such that
\begin{align*}
&\beta_0|l,m_1,m_2>=2l|l,m_1,m_2>, \\ 
&\alpha_0|l,m_1,m_2>=2m_1|l,m_1,m_2>,  \\
&\bar \alpha_0|l,m_1,m_2>=-2m_2|l,m_1,m_2>, \\
&N_ +|l,m_1,m_2>=0, {\text {and}} \\
&N_-|l,m_1,m_2>{\text { is a free $N_-$-module of rank 1.}}
\end{align*}

For each triple of complex numbers $r$, $s_1$ and $s_2$, we define the operator
$e^{2r\beta  +2s_1\alpha  +2s_2\bar \alpha }$ by the mapping 
$|l,m_1,m_2>$ to $ |l +r,m_1
+s_1,m_2 +s_2>$ such that it commutes with the action of $N_\pm$.
The normal ordering $: * :$ is defined according to 
$\alpha <\alpha_0$, $\bar \alpha<\bar
\alpha_0$, $\beta  <\beta_0$ and $N_-<N_ +$.

Consider the operators $X^\pm(z)\,: F_{l,m_1,m_2}\rightarrow
F_{l,m_1\pm1,m_2\pm1}$ defined by
\begin{align*}
X^{+}(z)&={1 \over{(q-q^{-1})}}:Y^+(z)\biggl\{Z_ +(q^{-{{k +2} \over2}}z)W_
+(q^{-{k \over2}}z)-W_-(q^{{k \over2}}z)Z_-(q^{{{k +2} \over2}}z)\biggr\}: \\
& ={1 \over{(q-q^{-1})}}(X_1^{+}(z)- X_2^{+}(z) ), \\
X^{-}(z)&={-1 \over{(q-q^{-1})}}:Y^-(z)\biggl\{Z_ +(q^{{{k +2} \over2}}z)W_
+(q^{{k \over2}}z)^{-1}-W_-(q^{-{k \over2}}z)^{-1}Z_-(q^{-{{k +2}
\over2}}z)\biggr\}:  \\
 &={1 \over{(q-q^{-1})}}(X_1^{-}(z)- X_2^{-}(z) ), \end{align*}
where

\begin{align*}
Y^{+}(z)
 &=\exp\left\{\sum_{m=1}^\infty q^{-{{km} \over 2}}{{z^m}
\over{[km]}}\bigl(\alpha_{-m} +\bar \alpha_{-m}\bigr)\right\}\\
 &\qquad e^{2(\alpha  +\bar \alpha)}z^{{1 \over k}(\alpha_0 +\bar
\alpha_0)}\exp\left\{-\sum_{m=1}^\infty q^{-{{km} \over2}}{{z^{-m}}
\over{[km]}}\bigl(\alpha_m +\bar \alpha_m\bigr)\right\},\\
Y^{-}(z)
 &=\exp\left\{-\sum_{m=1}^\infty q^{{{km} \over 2}}{{z^m}
\over{[km]}}\bigl(\alpha_{-m} +\bar \alpha_{-m}\bigr)\right\}\\
 &\qquad e^{-2(\alpha  +\bar \alpha )}z^{-{1 \over k}(\alpha_0 +\bar
\alpha_0)}\exp\left\{\sum_{m=1}^\infty q^{{{km} \over2}}{{z^{-m}}
\over{[km]}}\bigl(\alpha_m +\bar \alpha_m\bigr)\right\},
\end{align*}

\begin{align*}
Z_ +(z)&=\exp\left\{-(q-q^{-1})\sum_{m=1}^\infty z^{-m}{{[m]} \over{[2m]}}\bar
\alpha_{2}\right\}q^{-{1 \over 2}\bar \alpha_{0}},\\
Z_-(z)&=\exp\left\{(q-q^{-1})\sum_{m=1}^\infty z^{2}{{[m]} \over{[2m]}}\bar
\alpha_{-m}\right\}q^{{1 \over 2}\bar \alpha_ {0}}, 
\end{align*}

\begin{align*}
W_{+}(z)&=\exp\left\{-(q-q^{-1})\sum_{m=1}^\infty z^{-m}{{[m]}
\over{[2m]}}\beta_{2}\right\}q^{-{1 \over 2}\beta_{0}},\\  
W_{-}(z)&=\exp\left\{(q-q^{-1})\sum_{m=1}^\infty z^{2}{{[m]}
\over{[2m]}}\beta_{-m}\right\}q^{{1 \over 2}\beta_{0}}. 
\end{align*}

The action of those operators are well defined as explained in \cite{2} 
which give  us a bosonization of $U_q(\hat {\frak sl}(2))$ of level k.

First we  define two new operators: 

\begin{align*}
\bar Y^{+}(z)
 &=\exp\left\{\sum_{m=1}^\infty q^{-{{km} \over 2}}{{z^m}
\over{[km]}}\bigl(\bar \alpha_{-m}\bigr)\right\}\\
 &\qquad e^{2(\bar \alpha)}z^{{1 \over k}(\bar
\alpha_0)}\exp\left\{-\sum_{m=1}^\infty q^{-{{km} \over2}}{{z^{-m}}
\over{[km]}}\bigl(\bar \alpha_m\bigr)\right\},\\
\bar Y^{-}(z)
 &=\exp\left\{-\sum_{m=1}^\infty q^{{{km} \over 2}}{{z^m}
\over{[km]}}\bigl(\bar \alpha_{-m}\bigr)\right\}\\
 &\qquad e^{-2(  \bar \alpha )}z^{-{1 \over k}(\bar
\alpha_0)}\exp\left\{\sum_{m=1}^\infty q^{{{km} \over2}}{{z^{-m}}
\over{[km]}}\bigl(\bar \alpha_m\bigr)\right\},
\end{align*}

Let 
\begin{align*}
V^{+}(z)
 &=\exp\left\{\sum_{m=1}^\infty q^{-{{km} \over 2}}{{z^m}
\over{[km]}}\bigl(\alpha_{-m} \bigr)\right\}\\
 &\qquad e^{2(\alpha  )}z^{{1 \over k}(\alpha_0 
}\exp\left\{-\sum_{m=1}^\infty q^{-{{km} \over2}}{{z^{-m}}
\over{[km]}}\bigl(\alpha_m \bigr)\right\},\\
V^{-}(z)
 &=\exp\left\{-\sum_{m=1}^\infty q^{{{km} \over 2}}{{z^m}
\over{[km]}}\bigl(\alpha_{-m} \bigr)\right\}\\
 &\qquad e^{-2(\alpha   )}z^{-{1 \over k}(\alpha_0 
)}\exp\left\{\sum_{m=1}^\infty q^{{{km} \over2}}{{z^{-m}}
\over{[km]}}\bigl(\alpha_m \bigr)\right\},
\end{align*}
 
 We define two current operators $\Phi^{\pm}(z)$ as: 
\begin{defn}
 \begin{align*}
\Phi^{+}(z)
&={1 \over{(q-q^{-1})}}:\bar Y^+(z)\biggl\{Z_ +(q^{-{{k +2} \over2}}z)W_
+(q^{-{k \over2}}z)-W_-(q^{{k \over2}}z)Z_-(q^{{{k +2} \over2}}z)\biggr\}:, \\
\Phi^{-}(z)
&={1 \over{(q-q^{-1})}}:\bar Y^-(z)\biggl\{Z_ +(q^{{{k +2} \over2}}z)W_
+(q^{{k \over2}}z)^{-1}-W_-(q^{-{k \over2}}z)^{-1}Z_-(q^{-{{k +2}
\over2}}z)\biggr\}:.  \end{align*}
\end{defn}

We, then, have: 
\begin{gather*}
X^+(z)=V^+(z)\Phi^+(z),  \\
X^-(z)=V^-(z)\Phi^-(z). 
\end{gather*}
Here $V^{\pm}(z)$ commutes with  $\Phi^{\pm}(z)$. 
Then we know that the  two operators $\Phi^{\pm}(z)$
give 
us bosonic realization  of the   quantum 
parafermions associated to 
$U_q(\hat {\frak sl}(2))$ \cite{3}\cite{4} .

\section{ Main results}

We write down the correlation functions between the vertex operators. 
\begin{align*}
\bar Y^{+}(z)\bar Y^{-}(w) & = : \bar Y^{+}(z)\bar Y^{-}(w):
\exp \left\{ \Sigma {-(w/z)^m [2m] \over [km] }\right\} z^{4/k}
\\
&=: \bar Y^{+}(z)\bar Y^{-}(w):
\exp \left\{ \Sigma {(w/z)^m} {\frac 1 m} {- q^{(k-2)m}+q^{(k+2)m}
\over 1-q^{2km}} \right\} z^{4/k},
  \\
\bar Y^{-}(z)\bar Y^{+}(w) & = 
: \bar Y^{-}(z)\bar Y^{+}(w):
\exp \left\{ \Sigma {-(w/z)^m [2m] \over [km] }\right\} z^{4/k} \\
&=
: \bar Y^{-}(z)\bar Y^{+}(w): 
\exp \left\{ \Sigma {(w/z)^m} {\frac 1 m} {-q^{(k-2)m}+q^{(k+2)m}
\over 1-q^{2km}} \right\}z^{4/k} .
\end{align*} 
\begin{align*}
\bar Y^{+}(z) Z^{+}(w)& = : \bar Y^{+}(z)Z^+(w): 
,  \\
Z^+(z)\bar Y^{+}(w) & = 
: Z^+(z)\bar Y^{+}(w): \exp \left\{ {\Sigma (w/z)^m (q-q^{-1})
q^{-km \over 2} [m]\over m}
 \right\} z^{2}
\\&=: Z^+(z)\bar Y^{+}(w):
\exp \left\{ {\Sigma (w/z)^m (q^{-(km-2m)\over 2}-q^{-(km+2m)\over 2})
 \over m} \right\} z^{2} .
\end{align*} 
\begin{align*}
Z^-(z)\bar Y^{+}(w)& = :Z^-(z)\bar Y^{+}(w):
,  \\
\bar Y^{+}(z) Z^-(w)& = 
: Z^-(w)\bar Y^{+}(z): 
\exp \left\{ {\Sigma (w/z)^m (q-q^{-1})
q^{-km\over 2} [m]\over m}
 \right\} z^{2}
\\
&=: Z^-(w)\bar Y^{+}(z): 
\exp \left\{ {\Sigma (w/z)^m (q^{-(km-2m)\over 2}-q^{-(km+2m)\over 2})
\over m }\right\} z^{2} .
\end{align*} 
\begin{align*}
\bar Y^{-}(z) Z^{+}(w)& = : \bar Y^{+}(z)Z^+(w): 
,  \\
Z^+(z)\bar Y^{-}(w) & = 
: Z^+(z)\bar Y^{-}(w): 
\exp \left\{ {-\Sigma (w/z)^m(q-q^{-1})
q^{km\over 2} [m]\over m}
\right\} z^{2} \\
&= : Z^+(z)\bar Y^{-}(w):
\exp \left\{ { \Sigma (w/z)^m( q^{{km-2m}\over 2}-q^{{km+2m}\over 2}) \over m}
\right\} z^{2} .
\end{align*} 
\begin{align*}
Z^-(z)\bar Y^{-}(w) & = : Z^-(z)\bar Y^{-}(w):
,  \\
\bar Y^{-}(z) Z^-(w) & = 
: \bar Y^{-}(z) Z^-(w):
 \exp \left\{ {- \Sigma (w/z)^m(q-q^{-1})
q^{km \over 2}[m] \over m}
\right\} z^{2} \\
&=: \bar Y^{-}(z) Z^-(w):
\exp \left\{ { \Sigma (w/z)^m( q^{{km-2m}\over 2}-q^{{km+2m}\over 2}
\over m }\right\} z^{2} .
\end{align*} 
\begin{align*}
Z^-(z)Z^+(w) & = : Z^-(z)Z^+(w):
,  \\
Z^{+}(z) Z^-(w) & = 
: Z^{+}(z) Z^-(w) : 
 \exp \left\{ {\Sigma (w/z)^m (q-q^{-1})^2 [m]^2[km] \over [2m]m}
\right\} \\
&=: Z^{+}(z) Z^-(w) :  
\exp \left\{ {\Sigma (w/z)^m(-q^{4m}-1+2q^{2m})(q^{km}-q^{-km}) 
\over (1-q^{4m}) m}
\right\} 
\end{align*}
\begin{gather*}
W^-(z)W^+(w)  = : W^-(z)W^+(w):
,  \\
W^{+}(z) W^-(w)  = 
: W^{+}(z) W^-(w) :\exp \left\{ -{\Sigma (w/z)^m (q-q^{-1})^2 [m]^2 
[k+2m] \over [2m]m }\right\}  \\
\qquad = 
: W^{+}(z) W^-(w) : 
\exp \left\{ {\Sigma (w/z)^m (q^{4m}+1-2q^{2m})(q^{km+2m}-q^{-km-2m}) 
\over   (1-q^{4m}) m}
\right\} 
\end{gather*}
\begin{gather*}
\Phi^{-}(w)\Phi^{+}(z)=  
{1 \over{(q-q^{-1})}} {1 \over{(q-q^{-1})}} \\
\qquad\times
 [ :\bar Y^-(w)Z_ +(q^{{{k +2} \over2}}w)W_+(q^{{k \over2}}w)^{-1}: 
:\bar Y^+(z)Z_ +(q^{-{{k +2} \over2}}z)W_+(q^{-{k \over2}}z):
\\
\qquad\quad-:\bar Y^-(w)W_-(q^{-{k \over2}}w)^{-1}Z_-(q^{-{{k +2}\over2}}w):
:\bar Y^+(z)Z_ +(q^{-{{k +2} \over2}}z)W_+(q^{-{k \over2}}z):
\\
\qquad\quad- :\bar Y^-(w)Z_ +(q^{{{k +2} \over2}}w)W_+(q^{{k \over2}}w)^{-1}:
:Y^+(z)W_-(q^{{k \over2}}z)Z_-(q^{{{k +2} \over2}}z):
\\
\qquad\quad+:\bar Y^-(w)W_-(q^{-{k \over2}}w)^{-1}Z_-(q^{-{{k +2}\over2}}w):
:Y^+(z)W_-(q^{{k \over2}}z)Z_-(q^{{{k +2} \over2}}z): ]
\\ 
\qquad=\Sigma_{i=1,..4}V_i(w,z).
\end{gather*}
\begin{lemma} 
The correlation function of $V_i(w,z)$, $i\neq 4$ has a first 
order pole at $w=q^{k+2}z $, but the 
correlation function of  $V_4(w,z)$ has  neither zero nor a pole 
at  $w=q^{k+2}z $. 
\end{lemma} 

This follows from the correlation function formulas above. 
 
Let 
\[L(z)= \int_C \Phi^{-}(w)\Phi^{+}(z)dw, \]
where the contour C is around the point $w=zq^{k+2}$.

\begin{prop} The current operator $L(z)$ is  in the form: 
\[L(z)= \lambda_1(z)+\lambda_2(z)+\lambda_3(z), \] where 
$\lambda_i(z)$, for $i=1,2,3 $, is a vertex operator and {\small
 \begin{align*}
 \lambda_1(z)= & {1 \over{(q-q^{-1})^2}} 
:\bar Y^-(q^{k+2}z)Z_ +(q^{{{2k +6} \over2}}z)W_+(q^{{{3k+4} \over2}}z)^{-1} 
\bar Y^+(z)Z_ +(q^{-{{k +2} \over2}}z)W_+(q^{-{k \over2}}z): \\
\lambda_2(z)=& -{1 \over{(q-q^{-1})^2}}
:\bar Y^-(q^{k+2}z)W_-(q^{{{k+4} \over2}}w)^{-1}Z_-(q^{{{k +2}\over2}})
\bar Y^+(z)Z_ +(q^{-{{k +2} \over2}}z)W_+(q^{-{k \over2}}z): \\
\lambda_3(z)= & - {1 \over{(q-q^{-1})^2}} 
:\bar Y^-(q^{k+2}z)Z_ +(q^{{{2k +6} \over2}}z)
W_+(q^{{{3k+4} \over2}}z)^{-1}
Y^+(z)W_-(q^{{k \over2}}z)Z_-(q^{{{k +2} \over2}}z):.
\end{align*} 
}
\end{prop}

This follows from the lemma above. 

From this,  we know that $L(z)$ is a sum of three vertex operators.

Here we will present   two screening  operators basically introduced in \cite{2}. 
We define the operator $S^\pm(z)$ as:
\begin{gather*} 
S^+(z)=
 \exp\left\{\sum_{m=1}^\infty z^m{1 \over{[2m]}}\left(q^{{{k} \over
2}m}\beta_{-m} +q^{{{k +2} \over 2}m}\bar \alpha_{-m}\right)\right\} \\
 \qquad e^{(k +2)\bar\beta +k\bar \alpha} z^{{1 \over2}(\beta_{0} +\bar
\alpha_0)}\exp\left\{-\sum_{m=1}^\infty z^{-m}
{1 \over{[2m]}}\left(q^{{{k} \over
2}m}\beta_{2} +q^{{{k +2} \over 2}m}\bar \alpha_m\right)\right\}. 
\end{gather*}
\begin{gather*}
S^-(z)= 
\exp\left\{\sum_{m=1}^\infty z^m{1 \over{[2m]}}\left(q^{-{{k} \over
2}m}\beta_{-m}-q^{-{{k +2} \over 2}m}\bar \alpha_{-m}\right)\right\} \\
 \qquad e^{(k +2)\beta -k\bar \alpha} z^{{1 \over2}(\beta_{0} -\bar
\alpha_0)}\exp\left\{-\sum_{m=1}^\infty z^{-m}{1 \over{[2m]}}\left(q^{-{{k} \over
2}m}\beta_{2}-q^{-{{k +2} \over 2}m}\bar \alpha_m\right)\right\}. 
\end{gather*}
\renewcommand{\qed}{}

\begin{lemma}
 The correlation function of $S^\pm(z)S^\pm(z)$ are given as:
\begin{align*} 
S^+(z)S^+(w) &= : S^+(z)S^+(w): (z-w) , 
S^-(z)S^-(w) &= : S^-(z)S^-(w): (z-w). 
\end{align*}
$S^+(z)$ and $S^-(z)$ are fermions.
\end{lemma} 

We can also calculate first all the correlation functions of the 
products of these two operators with  component of  
$X^\pm(z)$, which are given as:
\begin{lemma}
\begin{align*}
S^+(z)X^+_1(w)& =  :S^+(z)X^+_1(w):(z-wq )^{-1} \\
S^+(z)X^-_1(w)&=   :S^+(z)X^-_1(w):(z-wq^{k+1}) \\
X^+_2(w)S^+(z)&=   :X^+_2(w)S^+(z): (w-zq)^{-1} \\
X^-_2(w)S^+(z)&=   :X^-_2(w)S^+(z): (w-zq^{k+1}) \\
S^-(z)X^+_1(w)& =  :S^-(z)X^+_1(w):(z-wq^{-k-1} ) \\
S^-(z)X^-_1(w)&=   :S^-(z)X^-_1(w):(z-wq^{-1})^{-1} \\
X^+_2(w)S^-(z)&=   :X^+_2(w)S^-(z): (w-zq^{-k-1}) \\
X^-_2(w)S^-(z)&=   :X^-_2(w)S^-(z): (w-zq^{-1})^{-1}
\end{align*}
\end{lemma} 

As shown in \cite{2}, we have 
\begin{lemma} 
\begin{align*}
&X^+(z)S^+(w)=-S^+(w)X^+(z)\sim {{\partial_q} \over{\partial_qw}} \left\{{1
\over {z-w}}Y^+(z) B_1^+(z)\right\}, \\
&X^-(z)S^+(w)=-S^+(w)X^-(z)\sim0,
\end{align*} 
where
\begin{align*}
B_1^+(z) 
 &=\exp\left\{\sum_{m=1}^\infty z^m{1 \over{[2m]}}\left(q^{{{k +2} \over
2}m}\beta_{-m} +q^{{{k +4} \over 2}m}\bar \alpha_{-m}\right)\right\} \\
 &\qquad e^{(k +2)\beta +k\bar \alpha} z^{{1 \over2}(\beta_{0} +\bar
\alpha_0)}\exp\left\{-\sum_{m=1}^\infty z^{-m}{1 \over{[2m]}}\left(q^{{{k +2} \over
2}m}\beta_{m} +q^{{{k +4} \over 2}m}\bar \alpha_m\right)\right\},
\end{align*}
and  
\[{{\partial_{p}} \over{\partial_{p}z}}f(z)={{f(p^{1/2}z)-f(p^{-1/2}z)} \over
{(p^{1/2}-p^{-1/2})z}},\]
for a function $f(z)$ on $\Bbb C^*$ and  a
scalar $p$ in $\Bbb C^*$. 
\end{lemma} 

Similarly we have
\begin{lemma} 
\begin{gather*}
X^+(z)S^-(w)=-S^+(w)X^+(z)\sim 0 \\
X^-(z)S^-(w)=-S^-(w)X^-(z) \\
\qquad\sim 
 {{\partial_q} \over{\partial_qw}} \left\{{1
\over {z-w}}: Y^-(z)Z_+(q^{(k+2\over 2})W_+(q^{k \over 
2}z) S^+(z):\right\}. 
\end{gather*}
\end{lemma}

\begin{prop} 
$\lambda_2(z)$ commutes with $S^-(z)$, 
\[\lambda_2(z)S^-(w)= : \lambda_2(z)S^-(w):, \]
and  $L(z)$ commutes with $S^-(z)$ up-to a total difference.  
$\lambda_3(z)$ commutes with $S^+(z)$, 
\[\lambda_3(z)S^+(w)= : \lambda_2(z)S^+(w):, \]
and  $L(z)$ commutes with $S^+(z)$ up-to a total difference.
\end{prop} 

For two current operators A(z) and B(z), by that  $A(z)$ commutes with $B(w)$ up-to a total difference we mean: 
\[ [A(z), b(w)]= r(\delta (\frac z w r_1)A(z)-\delta(\frac z w r_2) A(z)), \]
where  $r_1, r_2, r$ are  non-zero constants and $A(z)$ is a current operator.   

This follows from the lemmas above.

\begin{thm} 
The operator $L(z)$ gives a realization of the quantized  W-algebra of 
${\frak sl}(2,1)$. 
\end{thm} 

This follows from the fact that 
the generator for  the quantized W-algebra of 
${\frak sl}(2,1)$ is uniquely determined by the fact that it commutes
with two fermion up-to a total difference as explained  in the 
section above. Because $S^+(z)$ and $S^{-}(z)$ are both fermions, 
which  are not inverse of each other, therefore our $l(z)$ is exactly 
such an operator, which gives us a realization of 
 the quantized  W-algebra of 
${\frak sl}(2,1)$. 

In a subsequent paper, we plan to present the 
quantized W-algebra for ${\frak sl}(m,n)$. For some of those 
algebras, we expect that we can build similar constructions  using the 
operators , which  comes from the quotient of
the $\hat {\frak sl}(n+1) $ up to $\hat {\frak gl}(n)$. This 
should be related to the construction in \cite{8}. 
 
On the other hand, the W-algebra structure  constructed from the parafermions of 
$\hat {\frak sl}(2)$ was  studied from  a different point of 
 view\cite{10}\cite{11}, where the connection with $\hat {\frak sl}(2,1)$ was not utilized.  
 It is a very interesting question to look at the results of 
\cite{10}\cite{11} from the point view of the W-algebra of $ {\frak sl}(2,1)$.

\subsection*{Acknowledgement}
The authors would like to thank M. Jimbo  and A. Matsuo
for useful discussions.

\label{lastpage}

\end{document}